%Plain TeX file of dinitz.tex by Doron Zeilberger; Final version of
%Feb. 28, 1995. To appear in the Amer. Math. Monthly
%begin macros
\baselineskip=14pt
\parskip=10pt
\def\halmos{\hbox{\vrule height0.15cm width0.01cm\vbox{\hrule height
 0.01cm width0.2cm \vskip0.15cm \hrule height 0.01cm width0.2cm}\vrule
 height0.15cm width 0.01cm}}
\font\eightrm=cmr8  

\magnification=\magstephalf

\parindent=0pt
\overfullrule=0in
%end macros
\bf
\centerline
{The Method of Undetermined Generalization and Specialization}
\medskip\centerline
{Illustrated with}
\medskip\centerline
{ Fred Galvin's Amazing Proof of the Dinitz Conjecture} 
\rm
\bigskip
\centerline{ {\it Doron ZEILBERGER}\footnote{$^1$}
{\eightrm  \raggedright
Supported in part by the NSF. This paper will appear
in the Amer. Math. Monthly.
Version of Feb. 28, 1995.
First version: Dec. 21, 1994.
} 
}
 
At the very beginning of our waning century, in what turned out to be
the most influential mathematical address  ever delivered,
David Hilbert[H] said:
 
{\it ``If we do not succeed in solving a mathematical problem, the
reason frequently consists in our failure to recognize the more
general standpoint from which the problem before us appears only as 
a single link in  a chain of related problems.''}
 
One paragraph later, he also said:
 
{\it In dealing with mathematical problems, specialization plays, as
I believe, a still more important part than generalization.}
 
Alas, all this is easier said than done. How does
one find the `right' generalization and specialization? 
 
The answer is: just go ahead and start  proving the conjecture.
At first, leave the exact form of the generalization  and/or specialization
blank,  and as you go along, see what kind of generalization/specialization
would be required to make the proof work out. Keep  `guessing and   erasing'
until you get it done,  just like doing a  crossword puzzle.
 
I  will illustrate this proof strategy
in terms of Fred Galvin's[G] recent brilliant
proof of the Dinitz conjecture. Following a tradition that goes
back to  Euclid, Galvin presented his proof 
as a marvelous  but `static' completed  edifice, just like the solution to
yesterday's (or last Sunday's)
puzzle, that hides  all the trials and tribulations by which  it was arrived.
Not very useful for solving today's puzzle...
 
The Dinitz conjecture asserts that given $n^2$  arbitrary sets
$A_{i,j}$ ($1 \leq i,j \leq n$),  
each  having  $n$ elements, then it is always possible to
pick elements $a_{i,j} \in A_{i,j}$ such that $(a_{i,j})$ is a
`generalized  Latin square', which means that 
each  row and each column must have all its $n$ entries distinct.
 
In other  words,  given a party of  $n$ boys and  $n$  girls,
in which every boy must  dance  once  with every girl,  and
such that  each possible  couple $(i,j)$ knows how to dance (with each 
other\footnote{$^2$}
{\eightrm These are couples' dances and the two  dancers should  be
able to coordinate their steps, so it is possible  for Abe to be able
to dance the tango with Alice but not with  Barbara, although
Barbara may be able to dance it  with other boys.}) only 
$n$ dances, then  out of  the $n^{n^2}$ ways of  assigning
dances to couples, there is  at least one way in  which each of the $2n$ 
individuals dances a different dance in each of his or her $n$ performances.
 
Like many  people,  I first heard about
the Dinitz conjecture[ERT] when Jeannette Janssen[J] 
brilliantly   applied the powerful algebraic  method
of Alon and Tarsi[AT] to `almost' prove it:
she proved  the analogous  statement for `non-square'
rectangles.
 
As soon as I found  out about the Dinitz conjecture,  
I was  struck by  its simplicity.
Like so many  times before, it seemed to me that
there ought to be a  `simple' proof to such a  simple statement, and
I spent many hours  trying,  in vain,  to prove it.  
 
The reason I found my inability to prove the Dinitz conjecture  
so frustrating is that  it appears to be
`intuitively obvious'. When all the sets $A_{i,j}$ are (pairwise) disjoint,
then the  statement is obvious.  In  the other extreme,  when all
the sets $A_{i,j}$  coincide,  then  it   is also obvious: we have
the problem of constructing an  ordinary $n \times  n$ Latin  square. 
This can  be constructed by looking at
the  multiplication table of any group of order $n$, in particular, the
additive group of  $\{0,1, \dots ,n-1  \}$ mod $n$:
$$
\matrix{
0 &  1 &  \ldots & n-2 &n-1 \cr
1  & 2 & \ldots  & n-1 & 0 \cr
\vdots & \vdots & \vdots &  \vdots \cr
n-1 & 0 &  \ldots & n-3 & n-2  \cr} \quad .
\eqno(L)
$$
It is intuitively obvious that as there is less overlap between
the sets $A_{i,j}$, there would have to be  more options, which should make 
it even easier to construct a generalized Latin square.
 
The need to prove intuitively obvious statements is very  common
in analysis,  which is  why I find  it  such a forbidding subject.
Indeed,  this was one of the reasons  why, shortly after my
Ph.D., I decided to  change fields from analysis to combinatorics, which I
found  much more  gratifying, as there the gap between `convincing yourself'
and `convincing the referee' is usually so much smaller.
 
Last spring, while I was still spending an hour a day trying to prove
the Dinitz  conjecture,
I got an E-mail message from Herb Wilf, who forwarded  an  E-mail
message from Richard Ehrenborg, who forwarded  an  E-mail message
from  Jeannette Janssen, who forwarded an E-mail message from
Gil Kalai[K]\footnote{$^3$}
{\eightrm I  am pleased to be in the  same connected component as
Kalai, but I wish that Gil would draw a directed edge between
him   and myself, especially since there already exists an  edge
from me to him.}
which  contained a lucid and concise two-and-a-half page outline
of Fred Galvin's proof, that Kalai had compiled.
 
When I finished reading
and digesting the proof, I kicked myself.  I felt that I could have found
it myself,
had I only followed  Hilbert's advice, coupled with
the `crossword methodology' alluded  to  above. With the very generous
help of Lady Hindsight,
I will now describe
how I  (and you!)  could have, and  should have, found
the very same proof, {\it without any prior  knowledge of  combinatorics
or graph  theory.}
 
Hopefully, this  presentation would enable you (or,
better still, me!), to give an elementary
proof (not using 
`the class number  formula for the Selmer group associated to the 
symmetric square  representation of a  modular lifting')
of  Wiles' theorem,  and  an 
elementary (and `the first')  proof  of the Riemann Hypothesis.\footnote{$^4$}
{\eightrm 
RH is (almost) equivalent to the following elementary  statement: Let
$a_n$ be the difference between  the
number of square-free integers  between $1$ and $n$  with
an  even  number  of prime factors and the number of those  with an
odd number of prime  factors, then for  some constant $A$, 
${\eightrm a_n \leq   A n^{9999/10000}}$. 
This would already make you rich  and famous.
The full RH is equivalent  to replacing  the 9999/10000 by any
number  larger  than 1/2.
}
 
{\bf The method of undetermined coefficients}
 
Every mathematician (and electrical engineer) knows that in order
to find a particular solution of a linear differential equation such as
$$
y'' +y = x^2 \quad ,
$$
one writes $y=A x^2+Bx+C$, for some {\it undetermined} constants
$A,B,C$. We don't know yet what they are, but we hope that they
exist, and are constants. Assuming this, we plug it into the equation,
getting
$$
 A x^2 + B x + (C+2A)= x^2 \quad .
$$
 
Comparing the coefficients of $x^2 , x^1$ and $x^0$  on both sides,
leads to the system of equations $A=1, B=0, C+2A=0$, which leads
to the solution $A=1, B=0, C=-2$. Hence $y=x^2-2$ is a solution
of the given differential equation.
 
{\bf The method of undetermined parameters in proofs}
 
In many proofs in number theory and elsewhere (e.g. [I] pp. 27-28),
we take parameters, say $t_0$ and $\rho$, fiddle with them,
and only at the end commit ourselves to a relation between them
(for example, [I], p.28, $\rho=1/( \log t_0 +2)$), that produces
the desired effect.
 
{\bf How to Generalize Dinitz's Conjecture?}
 
There is something too narrow and `square' about the statement of
the  Dinitz conjecture. A natural generalization
that comes to mind is to arbitrary graphs. 
Calling the elements of the sets $A_{i,j}$ `colors',  the task of
the  Dinitz conjecture is  to color  each cell  $(i,j)$ by
one  of the colors of the set  of colors  $A_{i,j}$ that the cell
is  allowed to use, in such a way that no two cells sharing  the  same  row,
or the  same column, can be  colored by the same color.
This  immediately brings to  mind graph  coloring. The  $n \times n$
discrete square  is an undirected  graph having the $n^2$ vertices
$\{  (i,j)  \vert 1 \leq  i,j \leq  n \}$, and  each vertex
$(i,j)$ is connected  to the  $2(n-1)$ vertices $(i,j'),  j' \neq j$
and $(i',j),  i'  \neq i$. A natural generalization would have
the form
 
{\it If G is any graph in a  class X (that includes squares) then
whenever   each vertex $v$  is  assigned  a set of  colors $A_{v}$,
and the  cardinalities  of  the sets  $A_v$ satisfy condition
$Y(G)$,  then it is possible to properly color the vertices of  $G$
so  that the color  of each vertex $v$ is drawn from the set
$A_v$.}
 
(A coloring of  a   graph is {\it proper} if two vertices joined
by  an edge always receive different colors.)
 
For the time being,  both the class  $X$, and the condition $Y$, are  left
{\it blank}. All  we  need  is that the class X contains  the graphs of squares
and the condition  $Y(G)$  becomes `having cardinalities $\geq n$'
when $G$ happens to be the $n \times  n$ square.

{\bf  Occam's Razor  and Specialization}
 
Properly coloring a graph  means that  for every two  vertices
$x$ and $y$ that are connected by an edge, we  require that:
 
{\it The colors assigned to $x$ and $y$ must differ.}
 
This statement really embodies   two statements:
 
{\it The color of $x$ is different from the color of $y$  AND
The color of $y$ is different from the color of $x$.}
 
What a waste! Following Occam's advice, we can drop either one
of the two statements. This leads to the idea of directing
the edges of our $n  \times n$  square-graph and  to consider the set of 
{\it  directed graphs}.  This  class might be easier to handle, since it has
more structure. For  any given  graph of $e$ edges,  there are $2^e$
ways to make it  a  directed  graph.  A  proper coloring of
a  directed  graph is assigning colors to each vertex  such  that
whenever there is an edge from vertex  $x$  to vertex  $y$,  the
color  assigned to $x$ must differ from the color assigned to $y$.
So  in  order to  properly color an undirected
graph,  all you  need  is  to  be able to color a  single one
of its many possible  directed  versions. 
So now we  have one more free parameter at our disposal: the 
way  to `orient'  the  graph of  the  $n \times n$ square.
Let's call this  orientation  $Z$. The proposed generalization/specialization
is now:
 
{\it If $G$ is any directed graph in a  class $X$ (that includes the squares
with orientation  $Z$) then
whenever each vertex $v$ is assigned  a set of  colors $A_{v}$,
and the  cardinalities  of  the sets  $A_v$ satisfy condition
$Y(G)$,  then it is possible to properly color the vertices of  $G$
so that the color  of each vertex $v$ is drawn from the set
$A_v$.}
 
We would be done  if we could find {\it  some} class $X$, 
{\it some} orientation $Z$,
and {\it   some} condition $Y$, 
that would enable  a  proof, such that the squares with orientation
$Z$ belong to $X$, and:
 
$$
Y( n \times n  \quad  \hbox{square with orientation Z} )=
[ \vert A_{i,j} \vert \geq n, \hbox{ for all} \,\, 1 \leq i , j \leq n] 
\quad .
$$
 
Our  best bet would be
an inductive proof, since graphs are so amenable  to  induction.
Such a proof would presumably consist in a recursive algorithm
to color the vertices that would  involve, at each step, getting
rid of some of the vertices  and edges, as well as of
some of the colors, thus shrinking  the  graph, 
that must stay in our class $X$, whatever it is, and
shrinking the sets $A_v$, in such a way that
condition $Y$,  whatever it  is, still holds.
 
But first let's impose some natural restrictions on
the orientation $Z$ of the graph  of the square. One  of the
great
principles of mathematics (and life)  is {\it symmetry} and {\it balance}
(e.g.  balancing  the budget.) The  number of neighbors of each
of the  $n^2$ vertices of the  $n  \times n$ square is $2n-2$.
When  we stick arrows  in the edges, it makes sense  to do it
in such  a  way that at each vertex
there  would be as many outgoing  edges as incoming
edges. So let's impose, tentatively of course, the
following condition on the  still elusive orientation Z:
 
{\it The orientation Z of the $n \times  n$ square
should  be such that  in  the resulting directed graph,   every vertex
has outdegree  $n-1$.}
 
(The {\it outdegree} of a vertex is the  number of edges coming  out of it.)
 
Now it is time to think of condition $Y$. The larger the
cardinality of the set $A_v$, the more  options we have  to  color
the vertex  $v$. On the other hand the larger the outdegree  of $v$,
the more restrictions  we have. Since more freedom should go hand
in hand with more responsibility, it makes sense that the condition
$Y$  regarding the cardinality of the set $A_v$  should be related
to the outdegree of the vertex $v$. 
Since 
the color  of any  vertex $v$ should be different than all its
(outgoing)  neighbors, that might
happen to be all distinct, the number of `optional colors' at  $v$,
i.e.  the cardinality of $A_v$, should  be  at  least  one more than
the outdegree of  $v$. But wait a minute! In  our `symmetric orientation'
$Z$, the outdegrees are all $n-1$ and in the statements of Dinitz's
conjecture  all the  cardinalities of the sets $A_v$ are  $\geq n$, one
more than the  outdegree. This  leads us to conjecture that  the
condition $Y=Y(G)$ should be: $\vert  A_v \vert \geq outdegree(v)+1$.
 
Plugging this (tentative!) condition $Y$ into the `undetermined
generalization' of the Dinitz  conjecture,  we are  lead to
the following statement:
 
{\it If G is any directed graph in a  class X (that includes the squares
with orientation  Z), then
whenever each vertex $v$ is  assigned a set of colors $A_v$
of cardinality $> outdegree(v)$,  it  is always
possible  to  properly color the graph in such  a way that the color
of $v$ is drawn from $A_v$.}
 
It now remains to find the class $X$ that will make the proof work,
and then make sure that there  is  an orientation  $Z$ of  the  
$n \times n$ square
such  that  the outdegree of  every vertex is $n-1$, and that  belongs to $X$.
 
It  is easy to see that the class  of  {\it all}  directed graphs is
too big (why?). On  the  other extreme  the  empty class  $X$  obviously
(and vacuously) satisfies  the  theorem, but  no orientation of the
square can  ever  belong to it, of course.
 
Anyway, let's  leave the nature  of the class  $X$ blank for now,
and try and prove the `generalized' Dinitz statement.  
Pick one of  the colors in the union  of  the  $A_v$'s,
let's call  it `red'.  We would  like
to color `red' at  least one of the vertices  that are allowed to
be colored `red', remove these vertices and their incident edges,
thereby getting a smaller graph to which we would like to apply induction.
In order for the induction to work,
the smaller graph  $G'$ must still  belong to the class X and satisfy  
condition
$Y$ (that the corresponding sets $A'_v$ will have cardinality strictly larger
then the outdegree of $v$  for every vertex $v$ in the reduced   graph $G'$.)  
 
When we pick a subset of the vertices to be colored `red', this subset
should be {\it independent}, i.e. no pair of its members can be connected
by an edge, or else the coloring would not be proper.
 
These vertices, that were colored `red', were chosen amongst all those
vertices $v$ that had the `red' option, i.e. for which `red' $\in A_v$.
All the other ones that had
`red' as  one of their  options,  but were {\it not} colored `red',
now lose  that option.  For induction to
work, we  need that  the reduced graph should still satisfy condition $Y$,
which means that these vertices, which are still waiting their turn to be
colored, but just lost one of their options,  should also  lose one of
their  (outgoing) neighbors.
The only way that this could happen is for the `frustrated red' vertices
to have had at least one neighbor amongst the `departing reds'.
Then having colored the `realized red' vertices `red', and having removed
them, leaves us a graph in which each of the `frustrated red' vertices
gets compensated for their loss of the `red' option, by getting
rid of (at least) one of their annoying (outbound) neighbors.
 
So in order  for  the difference between  the cardinality  of the sets
$A_v$ and the outdegree  of  $v$  to   be  still  $\geq 1$, we need  that
out of all the vertices  that  have  `red'  as one of  their options,
it  is  possible to pick an {\it independent} 
subset of vertices that would exercise that
option, in such a way that all  the other  vertices, that had `red'
as one of their options before, but did
not  use this option, would have an  edge leading to one of those  vertices
that  did get colored `red'. If this is the case,
removing  the  vertices that were
just colored `red', and  the   edges 
adjacent to them, would  then yield a smaller graph $G'$ ,
that should still  belong  to $X$, with correspondingly
smaller sets  $A_v$  that  still satisfy condition   $Y$.
 
Since we don't know  beforehand  which of  the  vertices  would  have
`red' (or later,  `green'  or any other  color)  as  one of their
options,  and also want property $X$ to be
`hereditary' (with respect to induced subgraphs),
we  should `leave our options  open'
and require  that {\it  any} subset  $S$  of vertices
should have  this  property of there always being an {\it independent}  subset
$S'  \subset S$  such  that every vertex in $S  -  S'$ has  an edge
directed toward some  vertex  of $S'$.  
This  is exactly  the property $X$ that we
have  been  looking for,  and the proof that we already   have, works
with  that property  $X$. So   now  we  can  formulate:
 
{\bf Definition:} {\it A  directed  graph $G$  has property  $X$ if
for every   subset of   vertices $S$  there is an
independent  subset $S' \subset  S$
such  that every vertex  in $S-S'$ has an  edge  directed  toward  a  vertex
of  $S'$.}
 
We   have  just  proved:
 
{\bf  The  `Trivializing' Generalization:} {\it Let $G$  be a  directed
graph  having  property $X$
(defined above). If  every vertex $v$ is given a set  of colors
$A_v$ whose cardinality exceeds the outdegree of $v$,
then  it is  always possible
to properly color $G$  in such a  way 
that the  color of $v$ is picked  from $A_v$.}
 
But is this indeed  a  generalization of the  statement   of  the
Dinitz  conjecture?  We  still  need  to find an orientation $Z$
of the graph  of the $n \times n$ square
such that every vertex   has  outdegree  $n-1$, and
that has property $X$.
 
One  of  the many possible  ways of picking $Z$ is  by
picking  the following orientation. Looking at the Latin
square $(L)$ given above,
the horizontal (vertical)
edges are directed from the smaller (larger) entries to larger (smaller)
ones. In other words:
 
$$
(i,j) \rightarrow (i',j) \quad if   \quad
[ \, (i+j-2)  \quad mod \quad   n \, ]
> [ \, (i'+j-2) \quad mod \quad n \, ] \quad ,
$$
 
$$
(i,j) \rightarrow (i,j') \quad if   \quad
[ \, (i+j-2)  \quad mod \quad  n \,] < 
[ \, (i+j'-2) \quad mod \quad n \, ] \quad .
$$
 
To  prove   property   $X$, Galvin  invokes the famous
Gale-Shapley `Stable  Marriage' theorem  ([GS],[PTW]), 
with the rows  representing
men, the columns  representing   women, and an arrow
from  $(i,j)$  to $(i',j)$ meaning  that Ms.  $j$ prefers Mr.  $i'$ to
Mr. $i$ while an arrow  from $(i,j)$  to $(i,j')$ 
meaning  that Mr. $i$ prefers Ms.  $j'$ to Ms.  $j$. Having  property $X$
is easily seen  to be  equivalent to  the existence of  a  stable
marriage,  even if  some  of  the  relationships  are removed,
because of the laws of  the land.\footnote{$^5$}
{\eightrm For example forbidding $(i,j)$ where Mr.  $i$  is
a Cohen and Ms. $j$ is a divorc\'ee.}
In this more general   situation,
it is no longer guaranteed that  everybody  gets married, but 
those who do, do  so without fear of being scorned.
 
I  believe that even if  the Gale-Shapley algorithm and/or
theorem did not  exist, it would not   have been too  hard
to either discover  it from scratch,  or prove by other  means
(e.g. induction)
that   there is  some orientation  $Z$ 
(in particular  the  one given above), that  satisfies  property
$X$. We  invite  the reader  to do this right now! \halmos
 
{\bf Postscript:} The true story
is even more amazing, and I hope that Galvin would write up
the story that he told me after he received the
first draft of this paper. 
Since this is {\it his} story, not mine, I will
not give it away, except to quote Noga Alon who said:
`The moral of the (true) story of how Galvin found his proof
is not to follow Hilbert, but to follow a simpler
adage: {\it Know where to look things up'}. Of course, just like Hilbert's
advice, this is easier said than done, and it takes someone like
Galvin to use this so effectively. It is interesting to note that the
right $Y$ and $Z$ were already present in [AT] and [J].
 
Noga Alon has informed me that Kalai's two-and-a-half page
expos\'e was based on a one-page description that Alon has
sent Kalai, and that Alon wrote up based on Galvin's letter
to him. I wish to thank Mireille Bousquet-M\'elou, Fred Galvin, Bruno Salvy
and Herb Wilf for helpful remarks on an earlier version.
 
{\bf Note:} John  Noonan,  of Temple University, has  written
Maple programs  that  implement the algorithm in Galvin's proof
and the  Gale-Shapley algorithm. They are available by anonymous
{\tt  ftp} to {\tt  ftp.math.temple.edu}  in  directory
{\tt pub/noonan}, or via Mosaic to 
{\tt http://www.math.temple.edu/$\tilde{\quad}$ noonan}.
 
{\bf References}
 
[AT] N. Alon and  M. Tarsi, {\it  Coloring and orientations
in graphs }, Combinatorica {\bf 12} (1992), 125-134
\smallskip
[ERT] P. Erd\"os, A.L. Rubin, and H. Taylor, {\it Choosability in Graphs},
Congr. Numer. {\bf 26}(1980), 122-157.
\smallskip
[GS] D. Gale and L.S. Shapley, {\it College admissions and the stability of
marriage}, Amer. Math. Monthly {\bf 69}(1962), 9-15.
\smallskip
[G] Fred  Galvin, {\it The list chromatic index of a
bipartite multigraph}, J. Combin. Theory Ser. B. {\bf 63}(1995), 153-158.
\smallskip
[H] David Hilbert, {\it Mathematical Problems}, Lecture
delivered  before  the  International Congress  of   Mathematicians
at  Paris in 1900, Bull. Amer. Math. Soc. {\bf 8}(1902), 437-479.
[Translation by Mary Winston Newson  of the German original
that appeared in {\it  G\"ottinger Nachrichten}, 1900,pp.
253-297.] [Reprinted in: `Mathematical developments arising
from Hilbert Problems',  F.E.Browder,  Ed.,`Proceedings of
Symposia in Pure  Mathematics' v.XXVIII-Part1, Amer.  Math. Soc.,
Providence, 1976,  pp.  1-34.]
\smallskip
[I] A. E. Ingham,  {\it `The distribution of prime  numbers'},
Cambridge  University Press, Cambridge, 1932.
\smallskip
[J]  Jeannette C.M. Janssen,
{\it The Dinitz problem solved for rectangles},
Bull. Amer. Math. Soc. (N.S.) {\bf 29}(1993),  243-249.
\smallskip
[K] Gil  Kalai, {\it E-mail message to Jeannette Janssen} (and  possibly
to  others, but not to me),  dated 28 Jan. 1994 11:26 GMT.
\smallskip
[PTW] George Polya, Robert E. Tarjan, and Donald R. Woods,
{\it  ``Notes on introductory combinatorics''}, Birkhauser, Basel
and Cambridge, 1983.

Department of Mathematics, Temple University,
Philadelphia, PA 19122, USA. 
 
E-mail: {\tt zeilberg@math.temple.edu}
 
WWW: {\tt http://www.math.temple.edu/$\tilde{\quad}$zeilberg}
 
Anon. ftp: {\tt ftp.math.temple.edu}, directory {\tt /pub/zeilberg}.
\bye